\documentclass[12pt,twoside]{article}
\usepackage[latin1]{inputenc}
\usepackage{amsmath}
\usepackage{amssymb,amsfonts}
\usepackage{graphicx}
\usepackage{times,amssymb,amscd}
\usepackage[USenglish]{babel}
\usepackage{verbatim}
\usepackage{enumerate}

\newcommand{\bc}{\mathbf{c}}

\newcommand{\bx}{\mathbf{x}}

\newcommand{\bu}{\mathbf{u}}

\newcommand{\Bu}{\boldsymbol{u}}

\newcommand{\cD}{\mathcal{D}}

\newcommand{\cB}{\mathcal{B}}

\begin{document}
\pagestyle{myheadings}
\markboth{\centerline{J.~Szirmai}}
{Horoball packings in the hyperbolic $n$-space}
\title
{Horoball packings to the totally asymptotic regular simplex in the hyperbolic $n$-space
}

\author{\medbreak \medbreak {\normalsize{}} \\
\normalsize Jen\H{o}  Szirmai \\
\normalsize  Budapest University of Technology and Economics\\
\normalsize Institute of Mathematics, Department of Geometry \\
\normalsize H-1521 Budapest, Hungary \\
\normalsize Email: szirmai@math.bme.hu }
\date{\normalsize (\today)}
%%%%%%%%%%%%%%%%%%%%%%%%%%%%%%%%%%%%%%%%%%%%
%%%%%%%%%%%%%%%%%%%%%%%%%%%%%%%%%%%%%%%%%%%%
\maketitle
%%%%%%%%%%%%%%%%%%%%%%%%%%%%%%%%%%%%%%%%%%%%
\begin{abstract}

In \cite{Sz11} we have generalized the notion of the simplicial density function for horoballs in the extended hyperbolic space 
$\overline{\mathbf{H}}^n, ~(n \ge 2)$,
where we have allowed {\it congruent horoballs in different types} 
centered at the various vertices of a totally asymptotic tetrahedron.  
By this new aspect, in this paper we study the locally densest horoball packing arrangements and their densities with
respect to totally asymptotic regular tetrahedra in hyperbolic $n$-space $\overline{\mathbf{H}}^n$ extended with its 
absolute figure, where the ideal centers of horoballs give rise to vertices of a totally asymptotic regular tetrahedron.

We will prove that, in this sense, {\it the well known B\"or\"oczky density upper bound for "congruent horoball" packings of $\overline{\mathbf{H}}^n$
does not remain valid for $n\ge4$,} but these locally optimal ball arrangements do not have extensions to the whole $n$-dimensional hyperbolic space. 
Moreover, we determine an explicit formula for the density of the above locally optimal horoball packings, allowing horoballs in different types.    
\end{abstract}

%%%%%%%%%%%%%%%%%%%%%%%%%%%%%%%%%%%%%%%%%%%
%%%%%%%%%%%%%%%%%%%%%%%%%%%%%%%%%%%%%%%%%%%

\newtheorem{theorem}{Theorem}[section]
\newtheorem{corollary}[theorem]{Corollary}
\newtheorem{lemma}[theorem]{Lemma}
\newtheorem{exmple}[theorem]{Example}
\newtheorem{defn}[theorem]{Definition}
\newtheorem{rmrk}[theorem]{Remark}
\newtheorem{proposition}[theorem]{Proposition}
\newenvironment{definition}{\begin{defn}\normalfont}{\end{defn}}
\newenvironment{remark}{\begin{rmrk}\normalfont}{\end{rmrk}}
\newenvironment{example}{\begin{exmple}\normalfont}{\end{exmple}}
\newenvironment{acknowledgement}{Acknowledgement}

%%%%%%%%%%%%%%%%%%%%%%%%%%%%%%%%%%%%%%%%%%%%%%%%%%%%%%%%%%%%%%%%%%%%

%==================================================================%
%                             the main article                               %
%==================================================================%

%%%%%%%%%%%%%%%%%%%%%%%%%%%%%%%%%%%%%%%%%%%%%%%%%%%%%%%%%%%%%%%%%%%%

\section{Introduction}

We consider the horospheres and their bodies, the horoballs. A horoball packing $\cB=\{B\}$ of $\overline{\mathbf{H}}^n$ is an 
arrangement of non-overlapping horoballs ${B}$ in $\overline{\mathbf{H}}^n$. 
The notion of local density of the usual ball packing can be extended for horoball packings $\cB$ of $\overline{\mathbf{H}}^n$. Let
$B \in \cB$, and $P \in \overline{\mathbf{H}}^n$ an arbitrary point. Then, $\rho(P,B)$ is defined to be the length of the
unique perpendicular from $P$ to the horosphere $S_h$ bounding $B$, where again $\rho(P,B)$ 
is taken negative for $P \in B$. The Dirichlet--Voronoi cell (shortly D-V cell) $\cD(B)$ of $B$ in $\cB$ is defined to be
the convex body
\begin{equation}
\cD_h = D(\cB,B) := \{ P \in \mathbf{H}^n ~|~ \rho(P,B) \le \rho(P,B'), ~ \forall B' \in \cB \}. \tag{1.1}
\end{equation}
Since both, $B$ and $\cD$, are of infinite volume, the usual concept of local density has to be
modified. Let $Q \in \partial{\mathbf{H}}^n$ denote the base point (ideal center at the infinity) of $B$, and interpret $S$ as a Euclidean
$(n - 1)$-space. Let $B_{n-1}(R) \subset S_h$ be an $n-1$-ball with center $C \in S_h$. Then,
$Q \in \partial {\mathbf{H}^n}$ and $B_{n-1}(R)$ determine a convex cone $C_n(R) :=cone(B^Q_{n-1}(R)) \in \overline{\mathbf{H}}^n$ with
apex $Q$ consisting of all hyperbolic geodesics through $B_{n-1}(R)$ with limiting point $Q$. With
these preparations, the local density $\delta_n(B, \cB)$ of $B$ to $\cD$ is defined by 
\begin{equation}
\delta_n(\cB, B) :=\overline{\lim_{R \rightarrow \infty}} \frac{vol(B \cap C_n(R))} {vol(\cD \cap C_n(R))}, \tag{1.2}
\end{equation}
and this limes superior is independent of the choice of the center $C$ of $B_{n-1}(R)$. 

{\it In \cite{Sz11}  we have refined the notion of the ,,congruent" horoballs in a horoball packing to the horoballs of the "same type" 
because the horroballs are in general
congruent in the hyperbolic space $\overline{\mathbf{H}}^n$. 

Two horoballs in a horoball packing are in the "same type" if and only if the 
local densities of the horoballs to the corresponding cell (e.g. D-V cell; or ideal simplex, later on) are equal.} 
{\it \bf If we assume that the ,,horoballs belong to the same type"}, then by analytical continuation,
the well known simplicial density function on $\overline{\mathbf{H}}^n$ can be extended from $n$-balls of radius $r$ to the case $r = \infty$,
too. Namely, in this case consider $n + 1$ horoballs $B$ which are mutually tangent. The convex hull of their
base points at infinity will be a totally asymptotic or ideal regular simplex $T_{reg}^{\infty} \in \overline{\mathbf{H}}^n$ of finite
volume. Hence, in this case it is legitimate to write
\begin{equation}
d_n(\infty) = (n + 1)\frac{vol(B) \cap T_{reg}^\infty}{vol(T_{reg}^\infty)}. \tag{1.3}
\end{equation}
Then for a horoball packing $\cB$, there is an analogue of ball packing, namely (cf. \cite{B78}, Theorem
4)
\begin{equation}
\delta_n(\cB, B) \le d_n(\infty),~ \forall B \in \cB. \tag{1.4}
\end{equation}
\begin{rmrk}
The upper bound $d_n(\infty)$ $(n=2,3)$ is attained for a regular horoball packing, that is, a
packing by horoballs which are inscribed in the cells of a regular honeycomb of $\overline{\mathbf{H}}^n$. For
dimensions $n = 2$, there is only one such packing. It belongs to the regular tessellation $\{\infty, 3 \}$ . Its dual
$\{3,\infty\}$ is the regular tessellation by ideal triangles all of whose vertices are surrounded
by infinitely many triangles. This packing has in-circle density $d_2(\infty)=\frac{3}{\pi} \approx 0.95493.. $. 

In $\overline{\mathbf{H}}^3$ there is exactly one horoball packing whose Dirichlet--Voronoi cells give rise to a
regular honeycomb described by the Schl\"afli symbol $\{6,~3,~3\}$ . Its
dual $\{3,3,6\}$ consists of ideal regular simplices $T_{reg}^\infty$  with dihedral angle $\frac{\pi}{3}$ building up a 6-cycle around each edge 
of the tessellation.
\end{rmrk}

{\bf If horoballs of different types at the various ideal vertices are allowed}, then we can generalize the notion of the simplicial density function
\cite{Sz11}:
\begin{defn}
We consider an arbitrary totally asymptotic simplex $T=E_0 E_1$ $E_2 E_3 \dots E_n$ in the $n$-dimensional hyperbolic space $\overline{\mathbf{H}}^n$.
Centers of horoballs are required to lie at vertices of $T$. We allow horoballs $(B_i, ~ i=1,2,\dots,n)$ of different types at the various vertices and require to form a packing, 
moreover we assume that $$card(B_i \cap [E_{i_0}E_{i_1}\dots E_{i_{n-1}}]) \le 1, ~ ~ i_j \ne i, ~ ~ j \in \{0,1,\dots ,n-1\}.$$
(The hyperplane of points $E_{i_0},$ $E_{i_1},$ $\dots,$ $E_{i_{n-1}}$ is denoted by $[E_{i_0}E_{i_1}\dots E_{i_{n-1}}]$ may touch the horoball $B_{i_n}$.) 
The generalized simplicial density function for the above simplex and horoballs is defined as
\begin{equation}
\delta(\mathcal{B})=\frac{\sum_{i=0}^{n} vol(B_i \cap T)}{vol(T)}. \notag
\end{equation}
\end{defn}
In \cite{Sz11} I have studied the locally densest horoball packing arrangements and their densities with
respect to totally asymptotic tetrahedra $T(\alpha)$ in hyperbolic 3-space $\overline{\mathbf{H}}^3$, 
where the ideal centers of horoballs give rise to vertices of $T(\alpha)$.
Moreover, I have proved that, in this sense, {\it the well known B\"or\"oczky density upper bound for "congruent horoball" packings 
of $\overline{\mathbf{H}}^3$ does not remain valid.} 

In \cite{KSz} we have proved that the known B\"or\"oczky--Florian 
density upper bound for "congruent horoball" packings of $\overline{\mathbf{H}}^3$
remains valid for the class of fully asymptotic Coxeter tilings, even if
packing conditions are relaxed by allowing horoballs of different
types under prescribed symmetry groups. The consequences of this remarkable result are discussed for
various Coxeter tilings (see \cite{KSz}), and we have obtained four different optimal horoball packings with the maximal density.

\emph{Now, the main problem is} to find the  
locally densest horoball packing  related to the $n$-dimensional ($n\ge 4$) totally asymptotic {\it regular} simplex while allowing different
types of horoballs ($B_i$) to be centered at the vertices $E_i$ ($i=0,1,2,\dots n$) of the simplex, such that  
the density $\delta(\cB)$ (see Definition 1.2) of the corresponding horoball arrangement 
is maximal. In this case the horoball arrangement $\cB$ is 
said to be \emph{locally optimal.}

We will prove that, in this sense, {\it the well known B\"or\"oczky density upper bound for "congruent horoball" packings of $\overline{\mathbf{H}}^n$
$(n\ge 4)$
does not remain valid,} but these locally optimal ball arrangements do not have extensions to the whole $n$-dimensional hyperbolic space.  

For example, the density of this locally densest packing in $\overline{\mathbf{H}}^4$ is $\approx 0.77038$ which is larger than the B\"or\"oczky 
density upper bound $\approx 0.73046$.   

\section{Computations in projective model}

For $\overline{\mathbf{H}}^n$ $n \geq 2$ we use the projective model in Lorentz space
$\mathbf{E}^{1,n}$ of signature $(1,n)$, i.e.~$\mathbf{E}^{1,n}$ is
the real vector space $\mathbf{V}^{n+1}$ equipped with the bilinear
form of signature $(1,n)$
\begin{equation}
\langle ~ \mathbf{x},~\mathbf{y} \rangle = -x^0y^0+x^1y^1+ \dots + x^n y^n \tag{2.1}
\end{equation}
where the non-zero vectors
$$
\mathbf{x}=(x^0,x^1,\dots,x^n)\in\mathbf{V}^{n+1} \ \  \text{and} \ \ \mathbf{y}=(y^0,y^1,\dots,y^n)\in\mathbf{V}^{n+1},
$$
are determined up to real factors and they represent points in
$\mathcal{P}^n(\mathbf{R})$. ${\mathbf{H}}^n$ is represented as the
interior of the absolute quadratic form
\begin{equation}
Q=\{[\mathbf{x}]\in\mathcal{P}^n | \langle ~ \mathbf{x},~\mathbf{x} \rangle =0 \}=\partial \mathbf{H}^n \tag{2.2}
\end{equation}
in real projective space $\mathcal{P}^n(\mathbf{V}^{n+1},
\mbox{\boldmath$V$}\!_{n+1})$. All proper interior point $\mathbf{x} \in {\mathbf{H}}^n$ are characterized by
$\langle ~ \mathbf{x},~\mathbf{x} \rangle < 0$.

The points on the boundary $\partial \mathbf{H}^n $ in
$\mathcal{P}^n$ represent the absolute points at infinity of $\overline{\mathbf{H}}^n$.
Points $\mathbf{y}$ with $\langle ~ \mathbf{y},~\mathbf{y} \rangle >
0$ lie outside of $\overline{\mathbf{H}}^n$ and are called outer points
of $\mathbf{H}^n $. Let $X([\mathbf{x}]) \in \mathcal{P}^n$ a point;
$[\mathbf{y}] \in \mathcal{P}^n$ is said to be conjugate to
$[\mathbf{x}]$ relative to $Q$ when $\langle ~
\mathbf{x},~\mathbf{y} \rangle =0$. The set of all points conjugate
to $X([\mathbf{x}])$ form a projective polar hyperplane
\begin{equation}
pol(X):=\{[\mathbf{y}]\in\mathcal{P}^n | \langle ~ \mathbf{x},~\mathbf{y} \rangle =0 \}. \tag{2.3}
\end{equation}
Hence the bilinear form by (2.1) induces a bijection
(linear polarity $\mathbf{V}^{n+1} \rightarrow
\mbox{\boldmath$V$}\!_{n+1})$
from the points of $\mathcal{P}^n$
onto its hyperplanes.

Point $X [\bold{x}]$ and the hyperplane $\alpha
[\mbox{\boldmath$a$}]$ are called incident if the value of the
linear form $\mbox{\boldmath$a$}$ on the vector $\bold{x}$ is equal
to zero; i.e., $\bold{x}\mbox{\boldmath$a$}=0$ ($\mathbf{x} \in \
\mathbf{V}^{n+1} \setminus \{\mathbf{0}\}, \ \mbox{\boldmath$a$} \in
\mbox{\boldmath$V$}_{n+1} \setminus \{\mbox{\boldmath$0$}\}$).
Straight lines in $\mathcal{P}^n$ are characterized by the
2-subspaces of $\mathbf{V}^{n+1} \ \text{or $(n-1)$-spaces of} \
\mbox{\boldmath$V$}\!_{n+1}$ (see e.g. in \cite{M97}).

In this paper we set the sectional curvature of $\overline{\mathbf{H}}^n$,
$K=-k^2$, to be $k=1$. The distance $s$ of two proper points
$(\mathbf{x})$ and $(\mathbf{y})$ is calculated by the formula:
\begin{equation}
\cosh{{s}}=\frac{-\langle ~ \mathbf{x},~\mathbf{y} \rangle }{\sqrt{\langle ~ \mathbf{x},~\mathbf{x} \rangle
\langle ~ \mathbf{y},~\mathbf{y} \rangle }} . \tag{2.4}
\end{equation}
The foot point $Y(\bold y)$ of the perpendicular, dropped from the point $X(\bold x)$ 
on the plane $(u)$, has the following form:
\[
\bold y=\bold x -\frac{\langle \bold x, \bold u \rangle}
{\langle \bold u, \bold u \rangle} \bold u.
\tag{2.5}
\]

The length $l(x)$ of a horocycle arc to a chord segment $x$ is determined
by the classical formula due to {{J.~Bolyai}}:
\begin{equation}
l(x)=2 \sinh{\frac{x}{2}}.  \tag{2.6}
\end{equation}

The volume of the horoball sectors in the $n$-dimensional hyperbolic space $\overline{\mathbf{H}}^n$ can be calculated by the formula (2.7) 
which is the generalization of the classical formula of {{J.~Bolyai}}
to higher dimensions (see \cite{V}). 
If the volume of the polyhedron $A$
on the horosphere is $\mathcal{A}$, the volume determined by $A$ 
and the aggregate of axes drawn from $A$
is equal to 
\begin{equation}
V=\frac{1}{n-1}\mathcal{A}. \tag{2.7}
\end{equation}
%%%%%%%%%%%%%%%%%%%%%%%%%%%%%%%%%%%%%%%%%%%%%%%%%%%%%%%%%%%%%%%%%

\section{Horoball packings and their simplicial densities}
\subsection{Formula for the classical simplicial horoball density}
An $n$-simplex in $\overline{\mathbf{H}}^n$ is regular if its symmetry group operates transitively on the $k$-dimensional faces $(0 \le k \le n-1)$.
In this case, it has a unique barycenter (the fixed point of the symmetry group), and all of its edge lengths and dihedral
angles are of equal measure, respectively. Let $T_{reg}^{\infty}(2\alpha^n) \subset \overline{\mathbf{H}}^n$ denote a regular $n$-simplex of dihedral angles
$2\alpha^n \in [0, \pi]$. For hyperbolic simplicial $n$-volume, there are explicit formulas in terms of the dihedral angles only for $n\le6$. 
In case $n=2$ the area of an ideal triangle is $\pi$ and for the volume 
of an ideal regular simplex $Vol(T_{reg}^{\infty}(2\alpha_\infty^n))$ $(n \ge 3)$ in the $n$-dimensional hyperbolic space is due to J. Milnor \cite{Mi94}.  
\begin{theorem}[J.~Milnor]
Let $\beta=\frac{n+1}{2}$. Then, the ideal regular simplex volume $Vol(T_{reg}^{\infty})$ in the $n$-dimensional hyperbolic space $(n\ge 3)$, is given by
\begin{equation}
\begin{gathered}
Vol(T_{reg}^{\infty})=\sqrt{n} \sum_{k=0}^\infty \frac{\beta (\beta+1)\dots (\beta+k-1)}{(n+2k)!}A_{n,k} \ \ 
\text{where}\\
A_{n,k}=\sum_{i_0+\dots+i_n=k, \ \  i_k \ge 0}\frac{(2i_0)! \dots (2i_n)!}{i_0! \dots i_n!}. \tag{3.1}
\end{gathered}
\end{equation}
\end{theorem}
Th density of $n+1$ mutually tangent horoballs $B_i$ (the horoballs are in the same type) with respect to the regular simplex 
$Vol(T_{reg}^{\infty}(2\alpha_\infty^n))$ 
formed by their base points is given by (1.3). The classical universal density upper bound that is due to K.~B\"or\"oczky can be derived from 
this arrangemet in the the $n$-dimensional hyperbolic space $\overline{\mathbf{H}}^n$. In \cite{K98} R.~Kellerhals have proved the following 
\begin{theorem}[R.~Kellerhals]
The simplicial horoball density $d_n(\infty)$ $(n\ge3)$ is given by
\begin{equation}
d_n(\infty)=\frac{n+1}{n-1} \frac{n}{2^{n-1}}\prod_{k=2}^{n-1}\Big(\frac{k-1}{k+1}\Big)^{\frac{n-k}{2}} \frac{1}{Vol(T_{reg}^{\infty}(2\alpha_\infty^n))}
\tag{3.2}
\end{equation}
\end{theorem}
\begin{rmrk}
In the hyperbolic plane $\overline{\mathbf{H}}^2$ $d_2(\infty)=\frac{3}{\pi}$.
\end{rmrk}
\subsection{Formula for the generalized simplicial horoball density}
The aim of this section is to determine the optimal packing arrangements $\mathcal{B}_{opt}^n$ and their
densities for the totally asymptotic simplices in $\overline{\mathbf{H}}^n$ $(n \ge 2)$.
We will use the consequences of the following Lemma (see \cite{Sz05-1} in $3$-dimensions):

\begin{lemma}
Let $B_1$ and $B_2$ denote two horoballs with ideal centers $C_1$ and
$C_2$ respectively in the $n$-dimensional hyperbolic space $(n \ge 2)$. Take $\tau_1$ and $\tau_2$ to be two congruent
$n$-dimensional convex piramid-like regions, with vertices $C_1$ and $C_2$. Assume that these horoballs
$B_1(x)$ and $B_2(x)$ are tangent at point $I(x)\in {C_1C_2}$ and
${C_1C_2}$ is a common edge of $\tau_1$ and
$\tau_2$. We define the point of contact $I(0)$ such that the
following equality holds for the volumes of horoball sectors:
\begin{equation}
V(0):= 2 vol(B_1(0) \cap \tau_1) = 2 vol(B_2(0) \cap \tau_2). \notag
\label{szirmai-lemma}
\end{equation}
If $x$ denotes the hyperbolic distance between $I(0)$ and $I(x)$,
then the function
\begin{equation}
V(x):= vol(B_1(x) \cap \tau_1) + vol(B_2(x) \cap \tau_2)=\frac{V(0)}{2}(e^{(n-1)x}+e^{-(n-1)x}) \notag
\end{equation}
strictly increases as~$x\rightarrow\pm\infty$.
\end{lemma}

\textbf{Proof:} Let $\mathcal{L}$ and $\mathcal{L}'$ be parallel horocycles with centre $C$
and let $A$ and $B$ be two points on the curve $\mathcal{L}$ and $A':=CA \cap \mathcal{L}'$,
$B':=CB \cap \mathcal{L}'$.
By the classical formula of {J. Bolyai} 
$$
\frac{\mathcal{H}(A' B')} {\mathcal{H} (AB)}=e^{{x}}, 
$$
where the horocyclic distance between $A$ and $B$ is denoted by $\mathcal{H} (A,B)$. 

Then by the above formulas we obtain the following volume function:
$$
\aligned
 V(x)& =Vol (B_1 (x) \cap \tau_1)+Vol(B_2(x) \cap \tau_2)= \\
 & =\frac{1}{2} V(0)\Big(e^{(n-1)x}+\frac{1}{e^{(n-1)x}}\Big)
\endaligned 
$$
It is easy to see that this function strictly increases in the interval $(0,\infty)$. \ \ $\square$

We consider horoball packings with centers located at ideal vertices of an totally asymptotic regular simplex
$T_{reg}^{\infty}=E_0 E_1 E_2 \dots E_n$ in the $n$-dimensional hyperbolic space $\overline{\mathbf{H}}^n$ $(n \ge 2)$.

We consider the following two basic horoball configurations
$\mathcal{B}_n^i$, $(i=1,2)$:

\begin{enumerate}
\item All $n+1$ horoballs are of the same type and the adjacent horoballs
touch each other at the "midpoints" of each edge. This horoball arrangement is denoted by $\cB_0^n$. 
We define the point of tangency of two horoballs $B_0$ and $B_{n}$ on
side $E_0E_n$ to be $I(0)$ so that the following equality holds:
\begin{equation}
V(0):= (n+1) \cdot  Vol(B_0(0) \cap {T}_{reg}^\infty) = (n+1) \cdot Vol(B_{n}(0) \cap {T}_{reg}^{\infty})=(n+1) \cdot V_0. \notag
\end{equation}
In the hyperbolic plane ($n=2$) $V_0=1$ and we obtain applying the formula (3.2) that for dimensions $n\ge3$ 
$V_0=\frac{1}{n-1} \frac{n}{2^{n-1}}\prod_{k=2}^{n-1}\Big(\frac{k-1}{k+1}\Big)^{\frac{n-k}{2}}$ in the $n$-dimensional hyperbolic space. 
\item One horoball of the "maximally large" type centered at $E_n$.
The large horoball $B_n$ tangents the opposite face $E_0 E_1 E_2 \dots E_{n-1}$ of $T_{reg}^{\infty}$ and it determines the other $n$ horoballs 
touching the "large horoball". The point of tangency of $B_n$ and $B_0$ along
segment $I(0)E_0$ is denoted by $I(x_1)$ where $x_1$ is the hyperbolic distance between $I(0)$ and $I(x_1)$. This horoball arrangement is denoted by
$\cB_1^n$.
\end{enumerate}
\begin{figure}
\begin{center}
\includegraphics[width=6cm]{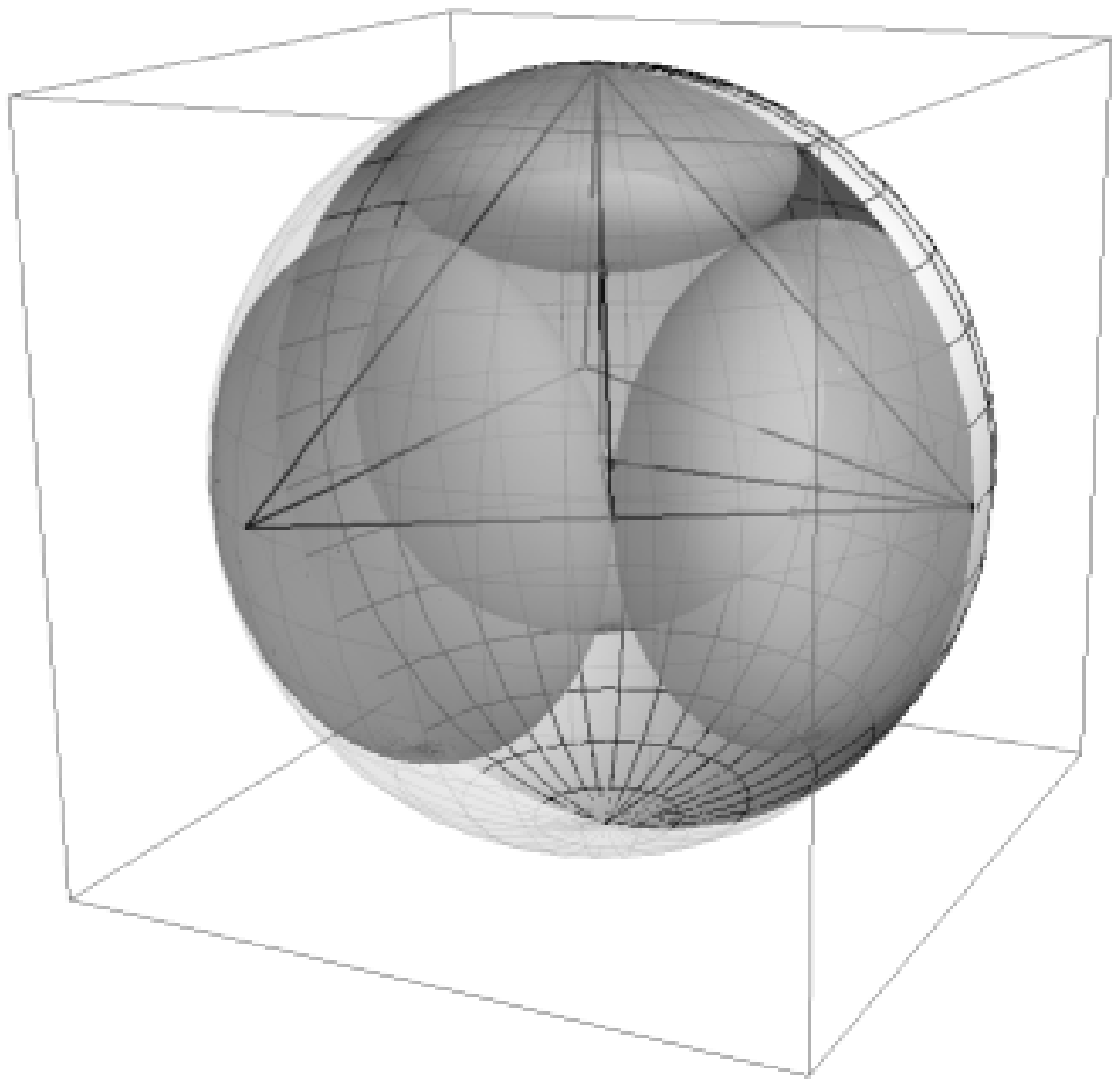}
\includegraphics[width=6cm]{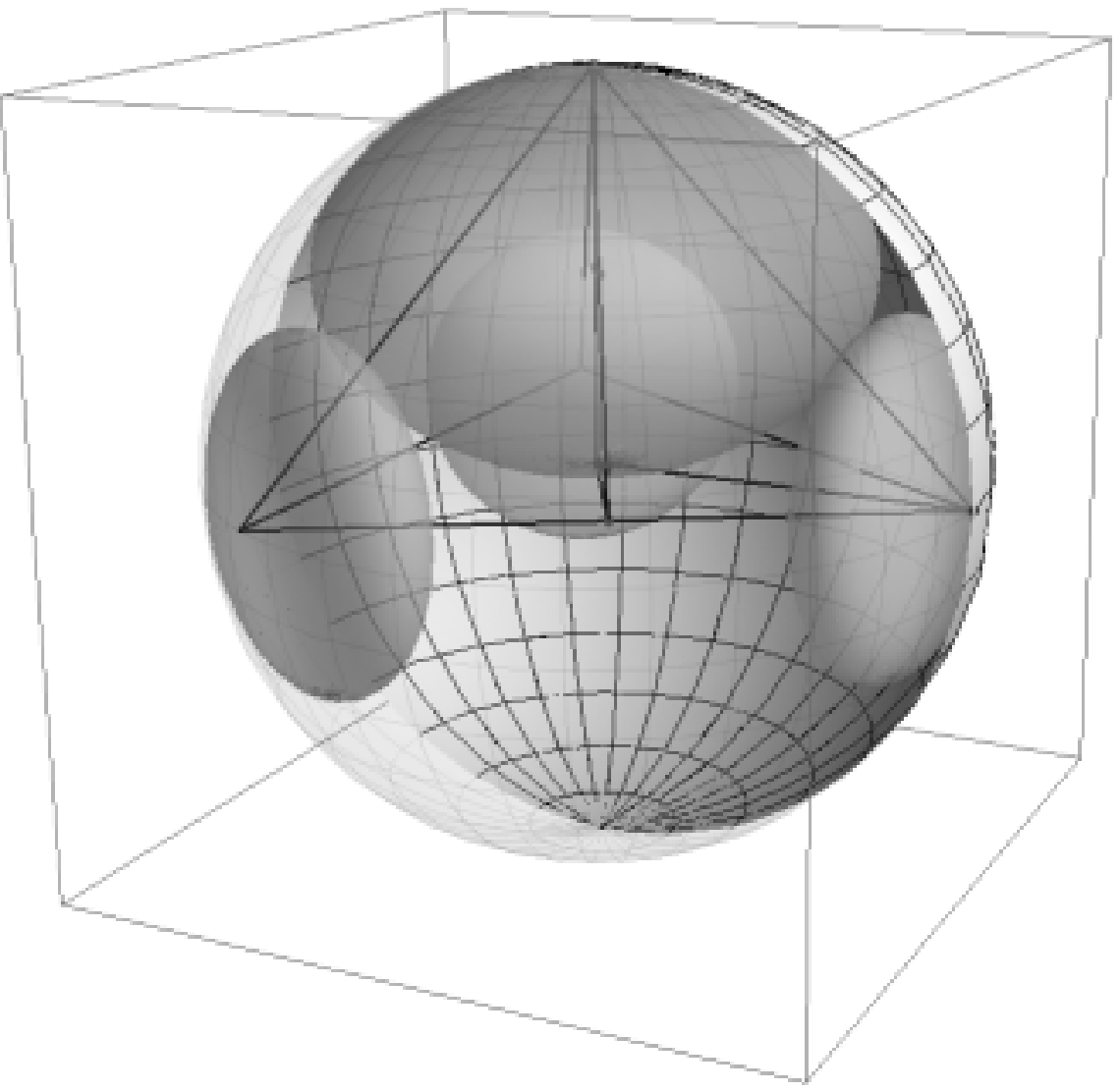}
\\a. ~~~~~~~~~~~~~~~~~~~~~~~~~~~~~~b. \\
\caption{Two optimal horoball arrangemets of $(3,3,6)$ tiling.}
\end{center}
\end{figure}
Due to symmetry considerations it is sufficient to consider the
cases when one horoball extends from the "midpoint" of an
edge until it touches the opposite side of the cell. 
Moreover, consider that point $I(x)$ is on edge $E_0E_n$.
This point is where the horoballs $B_i(x), \ (i=0,n)$ are tangent at
point $I(x) \in E_0I(0)$. Then $x$ is the hyperbolic distance
between $I(0)$ and $I(x)$. It is easy to see that we have to study that case where $x \in [0,x_1]$ and horoball $B_n$ 
touches the horoballs $B_i$ $(i=1,2,\dots, n-1)$.

In this case the function $V(x)$ can be computed by the
following formula
\begin{equation}
V(x):= n \cdot Vol(B_0(x) \cap {T}_{reg}^{\infty}) + Vol(B_3(x) \cap {T}_{reg}^{\infty}) \ \ x \in [0,x_1]. \notag
\end{equation}
Similarly to the Lemma 3.4, we can prove the following Lemma:
\begin{lemma}
\begin{equation}
\begin{gathered}
V(x):= n \cdot Vol(B_0(x) \cap {T}_{reg}^{\infty}) + Vol(B_3(x) \cap {T}_{reg}^{\infty})= \\
=V_0 (e^{(n-1) x} + n \cdot e^{-(n-1)x}), \ \ \ x \in [0,q_n], \notag
\end{gathered}
\end{equation}
and the maxima of function $V(x)$ are realized in point $I(0)$ if $q_n \le \frac{1}{n-1}\log{n}$ or at the point
$I(q_n)$ if $q_n \ge \frac{1}{n-1}\log{n}$ $(n\ge2)$.
\end{lemma}
\textbf{Proof:}
The second derivative of $V(x)$ is positive for all $n \ge 2$, thus it is strictly convex function and so, its maximum is achieved at 
$I(0)$ or at the point
$I(x_1)$. Moreover, $V(0)=V(\frac{1}{n-1}\log{n})$ and if $x=\frac{1}{2(n-1)}\log{n}$ then $V(x)$ is minimal.  \ $\square$

We consider a totally asymptotic regular tetrahedron ${T}_{reg}^{\infty}=E_0E_1 \dots E_n$ and place the horoball centers at vertices $E_0, \dots,
E_n$. We vary the types of the horoballs so that they satisfy our constraints of non-overlap.
The packing density is obtained by Definition 1.2. The dihedral angles of the above tetrahedron at the edges are equal of measure
$2\alpha_\infty^n=\arccos\Big(\frac{1}{n-1}\Big)$ (see \cite{K98}).

We introduce a Euclidean projective coordinate system to the tetrahedron ${T}_{reg}^{\infty}$ such that the center of the face 
$E_0 E_1 E_2 \dots, E_{n-1}$ coincide with the center of the model $O(1,0,0,\dots,0,0)$ moreover, we assume that $E_0 \sim (1,0,0,\dots,1,0)$
and $E_n \sim (1,0,0,\dots,0,1)$.

R. Kellerhals in \cite{K95} have proved, that the in-radius $\rho_n$ of an $n$ dimensional regular ideal simplex can be computed by the next formula:
\begin{equation}
\cosh(\rho_n)=\sqrt{\frac{2n}{n+1}}\cos(\alpha_\infty^n)=\frac{n}{\sqrt{(n-1)(n+1)}}, \tag{3.3}
\end{equation}
therefore the coordinates of the center \ $C \sim \alpha(1,0,0,\dots,0,0)+(1,0,0,\dots,0,1)=(\alpha+1,0,0,\dots,0,1)$ of ${T}_{reg}^{\infty}$ 
can be determined by formulas (2.4) and (3.3):
\begin{equation}
\begin{gathered}
\cosh{\rho_n}=\cosh{OC}=\frac{n}{\sqrt{(n-1)(n+1)}}=\frac{1}{\sqrt{1-\frac{1}{(\alpha+1)^2}}} \Rightarrow \\
C \sim (1,0,0,\dots,0,\frac{1}{n}). 
\end{gathered} \tag{3.4}
\end{equation}
\begin{figure}
\begin{center}
\includegraphics[width=10cm]{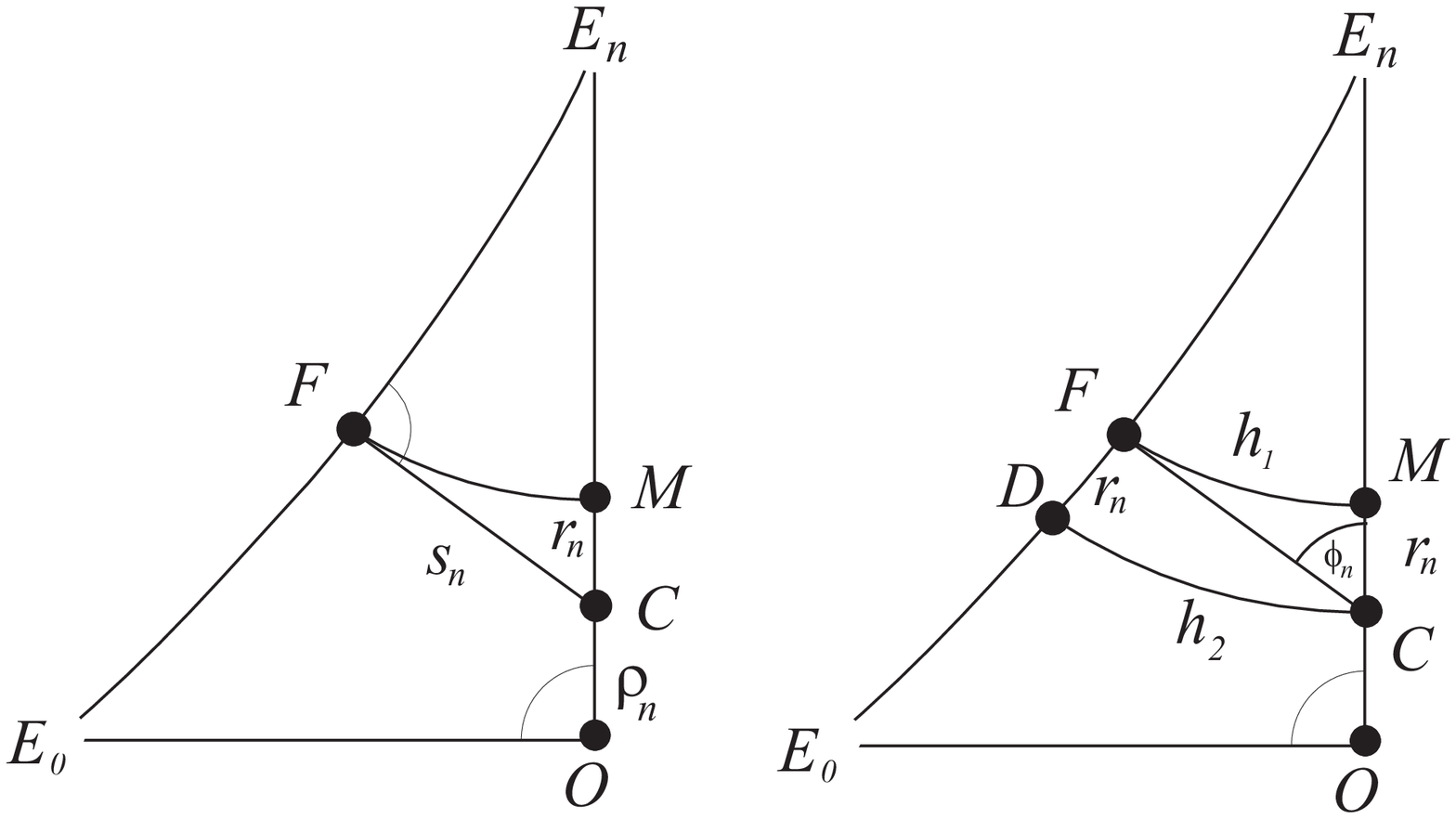}
\\a. ~~~~~~~~~~~~~~~~~~~~~~~~~~~~~~~~~~~~~~~~~b. \\
\caption{}
\end{center}
\end{figure}
Let the foot-point of perpendicular dropped from $C$ onto the staight line $E_0E_n$ be denote by $F$ which is the common point of the horoballs 
$B_n \in \cB_0$ and $B_0 \in \cB_0$ centered at $E_n$ and $E_0$. The hyperbolic distance $s_n=CF$ between the point $C(\bc)$ and the straight line
$E_0E_n$, given by $(\Bu)$, can be computed by the following formula (see Fig.~2):
\begin{equation}
\sinh(s_n)=\frac{|\langle \bx, \bu \rangle|}{\sqrt{-\langle \bx, \bx \rangle \langle \bu, \bu \rangle}}, \tag{3.5}
\end{equation}
where $\bu$ is the pol of the line $E_0E_n$ if we restrict the polarity to the plane $E_0E_nO$ (see Fig.~2).
The parallel distance of the angle $\phi_n=E_nCF\angle$ is $s_n$ therefore we obtain by the classical formula of J.~Bolyai and by formula (3.5) the 
following equation (see Fig.~2):
\begin{equation}
\sinh(s_n)=\cot(\phi_n)=\frac{n-1}{\sqrt{(n+1)(n-1)}}. \tag{3.6}
\end{equation}

We consider a horocycle $\mathcal{H}_2$ throught the point $C$ with center $E_n$ in the plane $E_nE_0O$ and the point $\mathcal{H}_2 \cap E_0E_n$ is denoted by $D$.
The horocyclic distances between points  $F$, $M$ and $C$, $D$ are denoted by $h_1$ and $h_2$. By means of formula of J.~Bolyai using the formula (3.6) 
yields
\begin{equation}
\begin{gathered}
\frac{h_2}{h_1}=e^{r_n}=\frac{1}{\sin(\phi)} \ \Rightarrow \\
r_n=\log\Big(\frac{1}{\sin(\phi_n)}\Big)=\log\Big(\sqrt{\cot^2(\phi_n)+1}\Big)=\log\Big(\sqrt{\frac{2n}{n+1}}\Big).
\end{gathered} \tag{3.7}
\end{equation}

The in-radius $\rho_n$ of an $n$ dimensional regular ideal simplex is by the formula (3.3)
\begin{equation}
\rho_n=\log\Big(\sqrt{\frac{n+1}{n-1}}\Big), \tag{3.8}
\end{equation}
therefore, the hyperbolic distance between the points $O$ and $M$ is
\begin{equation}
q_n=r_n+\rho_n=\log\Big(\sqrt{\frac{2n}{n-1}}\Big). \tag{3.9}
\end{equation}
From the Lemma 3.5 follows, that {\it the generalized simplicial horoball density} for the regular totally asymptotic tetrahedra is maximal at the $\cB_0$ horoball
configuration if $q_n \le \frac{1}{n-1}\log{n}$ or at the $\cB_1$ horoball arrangement if $q_n \ge \frac{1}{n-1}\log{n}$.

It is easy to see, that if $n=2,3$ then $q_n=\log\Big(\sqrt{\frac{2n}{n-1}}\Big)=\frac{1}{n-1}\log{n}$ and if $n\ge4$ then
$q_n>\frac{1}{n-1}\log{n}$.

Finally, we obtain by Lemma 3.5 and formula (3.9) the main results:
\begin{theorem}
For $n=2,3$ the maximal generalized simplicial horoball density to the regular totally asymptotic simplex in the $n$-dimensional hyperbolic space 
is realized at horoball 
arrangements $\cB_0^n$ and $\cB_1^n$, as well and if $n\ge4$ then the maximal density is attained at the horoball configuration $\cB_1^n$.
\end{theorem}
\begin{theorem}
The maximal generalized simplicial horoball density for the regular totally asymptotic simplex in the $n$-dimensional hyperbolic space 
can be determined by the following formula $(n\ge3)$:
\begin{equation}
\begin{gathered}
\delta(\cB^n_{opt})=\frac{1}{n-1} \frac{n}{2^{n-1}}\prod_{k=2}^{n-1}\Big(\frac{k-1}{k+1}\Big)^{\frac{n-k}{2}}\cdot \\ 
\cdot \Big( \sqrt{\frac{2n}{n-1}}^{ \ (n-1)}
+ n \cdot \sqrt{\frac{2n}{n-1}}^{\ (-n+1)} \Big) \frac{1}{Vol(T_{reg}^{\infty}(2\alpha_\infty^n))}
\tag{3.10}
\end{gathered}
\end{equation}
where $Vol(T_{reg}^{\infty}(2\alpha_\infty^n))$ denotes the ideal regular $n$-simplex volume.
\end{theorem}

\begin{rmrk}
\begin{enumerate}
\item The optimal horocycle packings in the 2-dimensional hyperbolic space have the following densities: $\delta(\cB_0^2)=\delta(\cB_1^2)
=\frac{3}{\pi} \approx 0.95493.. $
\item In the $4$-dimensional hyperbolic space the classical B\"or\"oczky's upper bound related to the horoball arrangement $\cB_0^4$ 
is $\approx 0.73046$, but
the new optimal generalized simplicial horoball density at ball configuration $\cB_1^4$ is $\approx 0.77038$.
\end{enumerate}
\end{rmrk}

The above results show, that the discussion of the densest horoball packings in the $n$-dimensionalen hyperbolic space with congruent horoballs 
in different types is not settled.
\vspace{3mm}

Optimal sphere packings in other homogeneous Thurston geometries
represent a class of open mathematical problems. For
these non-Euclidean geometries only very few results are known
\cite{Sz07-2}, \cite{Sz10-1}. Detailed studies are the objective of
ongoing research.

%%%%%%%%%%%%%%%%%%%%%%%%%%%%%%%%%%%%%%%%%%%%%%%%%%%%%%%%%%%%%%%%%%%%%%%%%%%%%%%%%%%%%%%%%%%%%

%============================================================================%
%                                references                                  %
%============================================================================%

\end{document}